\documentclass[12pt]{article}
\usepackage{pifont}

\usepackage{amsfonts}
\usepackage{latexsym}
\usepackage{amsmath}
\usepackage{amssymb}
\usepackage{color}

\newcommand{\CA}{\mathcal{A}_{0,\infty}}
\newcommand{\cB}{\mathcal{B}}
\newcommand{\CF}{\mathcal{F}}
\newcommand{\CP}{\mathcal{P}}
\newcommand{\BE}{\mathbb{E}}
\newcommand{\BF}{\mathbb{F}}
\newcommand{\CE}{\mathcal{E}}
\newcommand{\BR}{\mathbb{R}}

\newcommand{\BN}{\mathbb{N}}
\newcommand{\CN}{\mathcal{N}}

\newcommand{\E}{\hat{\mathbb{E}}}
\newcommand{\CH}{\mathcal{H}}
\newcommand{\CLL}{C_{l,lip}(\mathbb{R}^n)}

\newcommand{\CB}{C_{b,lip}(\mathbb{R})}
\newcommand{\FI}{\varphi}
\newcommand{\GX}{\underline{\sigma}^2}
\newcommand{\GS}{\bar{\sigma}^2}
\newcommand{\KB}{K_{\beta}}
\newcommand{\KS}{K_{\bar{\sigma}}}
\newcommand{\KSE}{K_{\bar{\sigma}}^{\varepsilon}}
\newcommand{\e}{\varepsilon}
\newcommand{\BB}{\bar{B}}

\newtheorem{theorem}{Theorem}[section]

\newtheorem{corol}{Corollary}[section]

\newtheorem{example}{Example}[section]

\newtheorem{lemma}{Lemma}[section]
\newtheorem{definition}{Definition}[section]

\title{An Invariance Principle of G-Brownian Motion for the Law of the Iterated Logarithm under G-expectation}

\date{}

\author{Panyu Wu\thanks{E-mail: wupanyu@mail.sdu.edu.cn} and Zengjing Chen\thanks{E-mail: zjchen@sdu.edu.cn
}\\Department of Mathematics,\\Shandong University, Jinan, China}

\begin{document}

\maketitle

\begin{abstract}
In this paper, we present the general invariance principle of
G-Brownian motion for the law of the iterated logarithm under
G-expectation, that is, for G-Brownian motion $(B_t)_{t\geq0}$ in
G-expectation space $(\Omega,L_G^1,\E)$, for any $n\geq3$, let
$$\zeta_n(t)=(2n\log\log n)^{-{1\over2}}B(nt),\ \ \ \forall t\in[0,1],$$
$$K_{\beta}:=\{x(\cdot):x\in C([0,1]),x(0)=0,\int_0^1|\dot{x(t)}|^2dt\leq\beta^2\},\ \beta\in\BR^+,$$
and if $\E[B_1^2]=\GS,\ -\E[-B_1^2]=\GX$, then\\
(I) the sequence $(\zeta_n)_{n\geq3}$ is relatively norm-compact quasi-surely,\\
(II) $v\{C(\zeta_n)\subseteq K_{\bar{\sigma}}\}=1,$\\
(III) $v\{C(\zeta_n)\supseteq K_{\underline{\sigma}}\}=1,$ \\
(IV) $\forall \beta\in[\underline{\sigma},\bar{\sigma}],
V\{C(\zeta_n)=K_{\beta}\}=1,$ \\
where $C(\zeta_n)$ denotes the cluster of sequence
$(\zeta_n)_{n=3}^{\infty}$, $(V,v)$ is the conjugate capacities
generated by G-expectation $\E[\cdot]$.
\par
And we also give some examples as applications.\\
\par  $\textit{Keywords:}$ Invariance Principle; Law of the Iterated Logarithm; Capacity; Sub-linear expectation;
G-Brownian Motion\\
\par
\emph{MR(2000):}\ \ \ 60F17, 60G50
\end{abstract}

\section{Introduction}
The classical law of the iterated logarithm (LIL for short) as fundamental
limit theorems in probability theory play an important role in the
development of probability theory and its applications. The original
statement of the LIL obtained by Khinchine (1924) \cite{kh} is to a
class of Bernoulli random variables. Kolmogorov (1929) \cite{ko} and
Hartman-Wintner (1941) \cite{har} extended Khinchine's result to
large classes of independent random variables. Strassen (1964)
\cite{strassen} extended Hartman-Wintner's result to large classes
of functional random variables, it is well known as the invariance
principle for LIL which provide an extremely powerful tool in
probability and statistical inference.
\par
Starting with Strassen\cite{strassen}, a wide of literature has
dealt with extension of the invariance priciple for the classical
IID conditions are so strong that limit the applications of the
invariance principle. So many papers have been to find weak
dependence or nonstationary conditions which are not only enough to
imply the invariance principle but also sufficiently general to be
satisfied in typical applications, see for example
\cite{ei},\cite{mc1},\cite{mc2},\cite{ma},\cite{wu}.
\par
On the other hand, the key in the proofs of the invariance principle
is the additivity of the probabilities and the expectations. In
practice, such additivity assumption is not feasible in many areas
of applications because the uncertainty phenomena can not be modeled
using additive probabilities or additive expectations. As an
alternative to the traditional probability\ expectation, capacities
or nonlinear probabilities\ expectations (for example Choquet
integral, g-expectation) have been studied in many fields such as
statistics, finance and economics.
\par
Recently, motivated by the risk measures, super-hedge pricing and
model uncertainty in finance, Peng \cite{p2006}-\cite{p2010}
initiated the notion of independently and identically distributed
(IID) random variables under sub-linear expectations. He also
introduced the notion of G-normal distribution and G-Brownian motion
as the counterpart of normal distribution and Brownian motion in linear case respectively. Under this
framework, he proved one law of large numbers (LLN for short) and the
central limit theorems (CLT for short) \cite{p2009}. As well, Chen proved the strong LLN
\cite{lln} and LIL \cite{lil} in this framework. G-expectation space
is the most important sub-linear expectation space introduced by
Peng \cite{p2006}, which take the role of Winener space in classical
probability. Now more and more people are interested in
G-expectation space or sub-linear expectation space.
\par
A natural question is the following: Can the classical invariance principle
for LIL be generalized under G-expectation space?
The purpose of this paper is
to investigate the invariance principle of G-Brownian motion for LIL
adapting the Peng's IID notion under G-expectations space, of course
we can not use the $\mathcal{X}^2$ distribution as in
\cite{strassen}.
\par
The remainder of this paper is organized as follows. In section 2,
we recall some notions and properties in the G-expectation space,
and prove some lemmas which will be useful in this paper. In section
3, we state and prove the main results of this paper, that is the
invariance principle of G-Brownian motion for LIL under G-expectation.
In section 4, we give some examples as applications of the new
invariance principle.
\section{Notations and Lemmas}
In this section, we introduce some basic notations and lemmas.
First, we shall recall briefly the notion of sub-linear expectation
and IID random variables initiated by Peng in \cite{p2009}. Let
$(\Omega,\CF)$ be a measurable space and let $\CH$ be a linear space
of random variables on $(\Omega,\CF)$.
\begin{definition}\cite{p2009}
A function $\BE:\CH\rightarrow \BR$ is called a sub-linear expectation, if it satisfies the following four properties: for all $X,Y\in\CH$\\
$(1)$ Monotonicity: $\BE[X]\geq\BE[Y]$ if $X\geq Y;$\\
$(2)$ Constant preserving: $\BE[c]=c,\ \forall c\in\BR;$\\
$(3)$ Sub-additivity: $\BE[X+Y]\leq \BE[X]+\BE[Y];$\\
$(4)$ Positive homogeneity: $\BE[\lambda X]=\lambda \BE[X],\ \forall
\lambda\geq 0.$
\par
The triple $(\Omega,\CH,\BE)$ is called a sub-linear expectation
space.
\end{definition}
Given a sub-linear expectation $\BE[\cdot]$, let us denote the
conjugate expectation $\CE[\cdot]$ of sub-linear expectation
$\BE[\cdot]$ by
$$\CE[X]:=-\BE[-X],\qquad X\in\CH.$$
\begin{definition}\cite{p2009}
\textbf{Independence and Identical distribution}\\
In a sub-linear expectation space $(\Omega,\CH,\BE)$, a random
vector $Y\in\CH^n$ is said to be independent from another random
vector $X\in\CH^m$ under $\E[\cdot]$ if
$$\BE[\FI(X,Y)]=\BE[\BE[\FI(x,Y)]_{x=X}],\ \forall \FI\in C_{l,lip}(\BR^{m+n}).$$
\par
Let $X_1$ and $X_2$ be two n-dimensional random variables in
sub-linear expectation spaces $(\Omega_1,\CH_1,\BE_1)$ and
$(\Omega_2,\CH_2,\BE_2)$ respectively. They are called identically
distributed, denoted by $X_1\stackrel d= X_2$, if
$$\BE_1[\FI(X_1)]=\BE_2[\FI(X_2)]\ ,\forall \FI\in\CLL.$$
\par
If $\bar{X}$ is independent from X and $\bar{X}\stackrel d= X$, then
$\bar{X}$ is said to be an independent copy of $X$.
\end{definition}
\begin{definition}\cite{p2009}
\textbf{G-normal distributed}\\ A random variable $X$ on
a sub-linear expectation space $(\Omega,\CH,\BE)$ is called G-normal
distributed, denoted by $X\sim \CN(0,[\underline{\sigma}^2,\GS])$,
if $$aX+b\bar{X}\stackrel d=\sqrt{a^2+b^2}X,\ \forall a,b\geq0,$$
where $\bar{X}$ is an independent copy of $X$, $\GS=\BE[X^2]$ and
$\underline{\sigma}^2=\CE[X^2].$ Here the letter G denotes the
function $G(\alpha):={1\over 2}\E[\alpha
X^2]={1\over2}(\GS\alpha^+-\GX\alpha^-):\BR\rightarrow\BR.$
\end{definition}
\begin{lemma}\label{le2}\cite{lil}
Suppose $\xi$ is distributed to G normal
$\CN(0,[\underline{\sigma}^2,\GS])$, where
$0<\underline{\sigma}\leq\bar{\sigma}<\infty.$ Let $\phi$ be a even
and bounded continuous positive function, then for any $b\in\BR$,
$$e^{-\frac{b^2}{2\GX}}\CE[\phi(\xi)]\leq\CE[\phi(\xi-b)].$$
\end{lemma}
\begin{definition}\cite{hu}
\textbf{G-Brownian motion}\\
Let
$G(\cdot):\BR\rightarrow\BR,G(\alpha)={1\over2}(\GS\alpha^+-\GX\alpha^-)$,
where $0\leq\underline{\sigma}\leq\bar{\sigma}<\infty.$ A stochastic
process $(B_t)_{t\geq0}$ in a sub-linear expectation space
$(\Omega,\CH,\BE)$
is called a G-Brownian motion if the following properties are satisfied:\\
(i) $B_0(\omega)=0;$\\
$(ii)$ For each $t,s\geq0$, the increment $B_{t+s}-B_t$ is
$\CN(0,[s\underline{\sigma}^2,s\GS])$-distributed and is independent
to $(B_{t_1},B_{t_2},\cdots,B_{t_n})$, for each $n\in\BN$ and $0\leq
t_1\leq t_2\leq\cdots\leq t_n\leq t.$
\end{definition}
$\bullet$ In the rest of this paper, we denote by
$\Omega=C_0(\BR^+)$ the space of all $\BR$-valued continuous
functions $(\omega_t)_{t\in\BR^+}$ with $\omega_0=0$, equipped with
the distance
$$\rho(\omega^1,\omega^2):=\sum_{i=1}^{\infty}2^{-i}[(\max_{t\in[0,i]}|\omega_t^1-\omega_t^2|)\wedge1].$$
$\bullet$ For every $\omega\in \Omega$, define the canonical process
by $B_t(\omega)=\omega_t,\ t\geq0$. The filtration generated by the
canonical process $(B_t)_{t\geq0}$ is defined by
$$\CF_t=\sigma\{B_s,\ 0\leq s\leq t\},\ \CF=\cup_{t\geq0}\CF_t.$$
$\bullet$ For each fixed $T\in[0,\infty),$ we set
$\Omega_T:=\{\omega_{\cdot\wedge T}:\omega\in\Omega\},$
$$L_{ip}(\Omega_T):=\{\FI(B_{t_1\wedge T},\cdots,B_{t_n\wedge T}):n\in \BN,t_1,\cdots,t_n\in[0,\infty),\FI\in\CLL\},$$
$$L_{ip}(\Omega):=\cup_{n=1}^{\infty}L_{ip}(\Omega_n).$$
$\bullet$ For any given monotonic increase and sub-linear function
$G(\cdot):\BR\rightarrow\BR$, there exist
$\bar{\sigma}\geq\underline{\sigma}\geq0$ such that
$G(\alpha)={1\over2}(\GS\alpha^+-\GX\alpha^-)$. We can construct
(see \cite{p2006},\cite{p2007}) a consistent sublinear expectation
called \textbf{G-expectation} $\E[\cdot]$ on $L_{ip}(\Omega)$, such that
$B_1$ is G-normally distributed under $\E[\cdot]$ and for each
$s,t\geq0$ and $t_1,\cdots,t_n\in[0,t]$ we have
$$\E[\FI(B_{t_1},\cdots,B_{t_n},B_{t+s}-B_t)]=\E[\psi(B_{t_1},\cdots,B_{t_n})],$$
where $\psi(x_1,\cdots,x_n)=\E[\FI(x_1,\cdots,x_n,\sqrt{s}B_1)].$
Under G-expectation $\E[\cdot]$, the canonical process
$(B_t)_{t\geq0}$ is a G-Brownian motion.
\par
We denote the completion of $L_{ip}(\Omega)$ under the norm $\| X\|
_p:=(\E[|X|^p])^{{1\over p}}$ by $L_G^p(\Omega),p\geq 1.$ And we
also denote the extension by $\E[\cdot]$. In the sequel, we consider the
G-Brownian motion means the canonical process $(B_t)_{t\geq0}$ under
the \textbf{G-expectation space} $(\Omega,L_G^1(\Omega),\E)$.
\par
Let $\CA$ denote the set of all
$[\underline{\sigma},\bar{\sigma}]$-valued, $\CF_t$-adapted
processes on the interval $[0,1]$. For each fixed $\theta\in\CA$,
set $P_{\theta}$, the law of the process
$(\int_0^t\theta_sdB_s)_{t\geq0}$ under the Wiener measure $P$. We
denote by $$\CP=\{P_{\theta}:\ \theta\in\CA\},$$ and define
$$V(A):=\sup_{\theta\in\CA}P_{\theta}(A),\ v(A):=\inf_{\theta\in\CA}P_{\theta}(A),\ \forall A\in\cB(\Omega).$$
It is easy to check that $$V(A)+v(A^c)=1,\ \forall
A\in\cB(\Omega),$$ where $A^c$ is the complement set of $A$. Through
this paper, we assume that $(V,v)$ is the conjugate capacities
generated by G-expectation $\E[\cdot]$. From Denis et al.\cite{denis}, we
know that $\CP$ is tight. For each $X\in L^0(\Omega)$ (the space of
all Borel measurable real functions on $\Omega$) such that
$E_{\theta}(X)$ exists for any $\theta\in\CA$, define the upper expectation
$$\bar{\BE}[X]:=\sup_{\theta\in\CA}E_{\theta}(X).$$
From Denis et al. \cite{denis}, for all $X\in L_G^1(\Omega)$, it
holds that $\bar{\BE}[X]=\E[X].$
\begin{definition}\cite{denis}
\textbf{quasi-surely}\\
A set $D$ is polar set if $V(D)=0$ and a property holds
``quasi-surely" (q.s. for short) if it holds outside a polar set.
\end{definition}
\begin{lemma}\label{le1}\cite{lil} \textbf{Borel-Cantelli lemma} \\
Let $\{A_n,n\geq1\}$ be a sequence of events in $\CF$ and $(V,v)$ be
a pair of capacities generated by
G-expectation $\E[\cdot]$.\\
$(1)$ If $\sum_{n=1}^{\infty}V(A_n)<\infty$, then $V(\cap_{n=1}^{\infty}\cup_{i=n}^{\infty}A_i)=0.$\\
$(2)$ Suppose that $\{A_n,n\geq1\}$ are pairwise independent with
respect to $V$, that is
$$V(\cap_{n=1}^{\infty} A_n^c)=\prod_{n=1}^{\infty}V(A_n^c).$$
If $\sum_{n=1}^{\infty}v(A_n)=\infty$, then
$v(\cap_{n=1}^{\infty}\cup_{i=n}^{\infty}A_i)=1.$
\end{lemma}
\begin{lemma}\label{le40}
$(1)$ Let $X_1$ and $X_2$ be two real random variables in
G-expectation spaces $(\Omega,L_G^1,\E)$. If $X_1\stackrel d= X_2$,
then
$$\bar{\BE}[\FI(X_1)]=\bar{\BE}[\FI(X_2)],\ i.e.\ V(X_1\in
A)=V(X_2\in A),$$ where $\FI(x)=I_A(x),\ A$ is an interval of the following types: $(-\infty,a)$, $(-\infty,a]$, $(a,b)$, $(a,b]$,
$[a,b)$, $[a,b]$, $[b,\infty)$, or $(b,\infty)$ with $(a,b)\in
\BR^2$ and $a<b$.
\par
$(2)$ In G-expectation space $(\Omega,L_G^1,\E)$, a real random
variable $Y$ is independent from another real random variable $X$
under $\E[\cdot]$, then
$$\bar{\BE}[\FI(X,Y)]=\bar{\BE}[\bar{\BE}[\FI(x,Y)]_{x=X}],\ i.e.\ V(X\in A,Y\in B)=V(X\in A)V(Y\in B),$$
where $\FI(x,y)=I_A(x)I_B(y)$, $A$ is an interval of the following
types: $(-\infty,a)$, $(-\infty,a]$, $(a,b)$, $(a,b]$, $[a,b)$,
$[a,b]$, $[b,\infty)$, or $(b,\infty)$; $B$ is an interval of the following types: $(-\infty,c)$, $(-\infty,c]$, $(c,d)$, $(c,d]$,
$[c,d)$, $[c,d]$, $[d,\infty)$ or $(d,\infty)$ with $(a,b,c,d)\in
\BR^4$ and $a<b,c<d$.
\end{lemma}
\textbf{Proof.} (1). We just consider the case $a=y,\
A=(-\infty,y]$, other cases can be proved similarly. Let
$\BF(y):=\bar{\BE}[I_{X_1\leq y}]=V(X_1\leq y)$, then $\BF$ is a
continuous function on $\BR$ from lemma 8 in page 143 of
\cite{denis}. So for each $y\in\BR$, for any $\e>0$, there exists
$\delta>0$ such that
$$|\BF(\bar{y})-\BF(y)|<\e,\ \forall \bar{y}\in[y-\delta,y+\delta].$$
Now we define two auxiliary functions
\begin{displaymath}
f(x)= \left\{\begin{array}{cc} 1 & x\in(-\infty, y-\delta);\\
\frac{(y-x)}\delta & x\in[y-\delta,y];\\
0 & x\in(y,\infty),
\end{array} \right.
g(x)= \left\{\begin{array}{cc} 1 & x\in(-\infty, y);\\
\frac{(y+\delta-x)}\delta & x\in[y,y+\delta];\\
0 & x\in(y+\delta,\infty).
\end{array} \right.
\end{displaymath}
Because $X_1\stackrel d= X_2$, using the monotonicity of
$\bar{\BE}[\cdot]$, we get
$$\BF(y-\delta)\leq\E[f(X_1)]=\E[f(X_2)]\leq\bar{\BE}[I_{X_2\leq y}]\leq \E[g(X_2)]=\E[g(X_1)]\leq \BF(y+\delta).$$
Therefore
$$-\e\leq\BF(y-\delta)-\BF(y)\leq \bar{\BE}[I_{X_2\leq y}]-\bar{\BE}[I_{X_1\leq y}]\leq\BF(y+\delta)-\BF(y)\leq\e.$$
By the arbitrariness of $\e$ we get $\bar{\BE}[I_{X_2\leq
y}]=\bar{\BE}[I_{X_1\leq y}]$, that is $V(X_1\leq y)=V(X_2\leq y)$.
\par
(2). We only consider the case $A=(-\infty,a],\ B=(-\infty,c]$,
other cases can be proved in the same way. For $\CP$ being compact,
we get $\BF(a,c):=\bar{\BE}[I_{X\in A}\cdot I_{Y\in B}]=V(X\leq
a,Y\leq c)$ is a continuous function on $\BR^2$ also from lemma 8 in page
143 of \cite{denis}. So for each $(a,c)\in\BR^2$, for any $\e>0$,
there exists $\delta>0$ such that
$$|\BF(\bar{a},\bar{c})-\BF(a,c)|<\e,\ \forall (\bar{a},\bar{c})\in[a-\delta,a+\delta]\times [c-\delta,c+\delta].$$
Similarly as (1) we define two auxiliary functions
\begin{displaymath}
f(x,y)= \left\{\begin{array}{cc} 1 & (x,y)\in(-\infty, a-\delta)\times (-\infty, c-\delta);\\
0 & (x,y)\in(a,\infty)\times (c,\infty);\\
\frac{(a-x)(c-y)}{\delta^2} & \hbox{others},
\end{array} \right.
\end{displaymath}
\begin{displaymath}
g(x,y)= \left\{\begin{array}{cc} 1 & (x,y)\in(-\infty, a)\times (-\infty, c);\\
0 & (x,y)\in(a+\delta,\infty)\times (c+\delta,\infty);\\
\frac{(a+\delta-x)(c+\delta-y)}{\delta^2} & \hbox{others}.
\end{array} \right.
\end{displaymath}
Obviously, $$f(x,y)\leq\FI(x,y)=I_{X\leq a}\cdot I_{Y\leq c}\leq
g(x,y.)$$ For $Y$ is independent from $X$ under $\E[\cdot]$, we get
$$\BF(a-\delta,c-\delta)\leq\E[f(X,Y)]=\E[\E[f(x,Y)]_{x=X}]\leq\bar{\BE}[\bar{\BE}[\FI(x,Y)]_{x=X}],$$
while
$$\bar{\BE}[\bar{\BE}[\FI(x,Y)]_{x=X}]\leq\E[\E[g(x,Y)]_{x=X}]=\E[g(X,Y)]\leq\BF(a+\delta,c+\delta).$$
It deduce that
$$-\e\leq\BF(a-\delta,c-\delta)-\BF(a,c)\leq \bar{\BE}[\bar{\BE}[\FI(x,Y)]_{x=X}]-\bar{\BE}[\FI(X,Y)]\leq\BF(a+\delta,c+\delta)-\BF(a,c)\leq\e.$$
By the arbitrariness of $\e$ we get $\bar{\BE}[\FI(X,Y)]=
\bar{\BE}[\bar{\BE}[\FI(x,Y)]_{x=X}]$, that is $V(X\leq a,Y\leq
c)=V(X\leq a)V(Y\leq c)$.\ $\square$
\par
\textbf{Remark:} A similar result as (1) of Lemma \ref{le40} can be
proved similarly when $X_1,X_2$ are in two sub-linear expectation
spaces respectively and one of the probability set generated by the
sub-linear expectation is weakly compact. Meanwhile, the conclusion
(2) of Lemma \ref{le40}, of course, holds true if $X,Y$ are in a
sub-linear expectation spaces and the probability set generated by
the sub-linear expectation is weakly compact, the proof is also
similarly.
\begin{lemma}\label{le3}
For any $s\leq t$, we have for almost surely $y\in\BR^+$
$$V(|B_s|\geq y)\leq V(|B_t|\geq y).$$
\end{lemma}
\textbf{Proof.} Let $(\Omega,\CF,P)$ be a probability space and
$(W_t)_{t\geq0}$ be Brownian motion in this space,
$(\CF_t)_{t\geq0}$ is the filtration generated by $(W_t)_{t\geq0}$.
Denis et al. \cite{denis} have proved
\begin{equation}\label{eq6}
\E[\FI(B_1)]=\sup_{\theta\in\Theta}E_P[\FI(\int_0^1\theta_sdW_s)], \
\forall \FI\in\CB,
\end{equation}
where $\Theta$ denote the collection of all
$[\underline{\sigma},\bar{\sigma}]$-valued $\CF_t$-adapted process
on interval [0,1]. To proof the lemma, we fist prove
\begin{equation}\label{eq7}
\bar{\BE}[I_{B_1\leq
y}]=\sup_{\theta\in\Theta}E_P[I_{\int_0^1\theta_sdW_s\leq y}].
\end{equation}
Define
$\BF(y):=\sup_{\theta\in\Theta}E_P[I_{\int_0^1\theta_sdW_s\leq y}],$
then $\BF(\cdot)$ is a not decreasing function on $\BR$. So
$\BF(\cdot)$ is almost surely continuous. Take $y$ is the continuous
point of $\BF(\cdot)$, for any $\e>0$, there exists $\delta>0$ such
that
$$|\BF(\bar{y})-\BF(y)|<\e,\ \forall \bar{y}\in[y-\delta,y+\delta].$$
Using the same auxiliary functions in the proof of the first part of
Lemma (\ref{le40}) and equality (\ref{eq6}), we get
$$\BF(y-\delta)\leq\E[f(B_1)]\leq\bar{\BE}[I_{B_1\leq y}]\leq \E[g(B_1)]\leq \BF(y+\delta).$$
Therefore
$$-\e\leq\BF(y-\delta)-\BF(y)\leq \bar{\BE}[I_{B_1\leq y}]-\BF(y)\leq\BF(y+\delta)-\BF(y)\leq\e.$$
By the arbitrariness of $\e$ we get $\bar{\BE}[I_{B_1\leq
y}]=\BF(y),$ that is equality (\ref{eq7}) hold.
\par
Similarly, we can show for almost surely $y\in\BR^+$ the following
equality (\ref{eq8}) hold
\begin{equation}\label{eq8}
\bar{\BE}[I_{|B_1|\geq
y}]=\sup_{\theta\in\Theta}E_P[I_{|\int_0^1\theta_sdW_s|\geq y}].
\end{equation}
Since $B_t\stackrel d=\sqrt t B_1$, from Lemma \ref{le40} and
equality (\ref{eq8}) we have
$$V(|B_t|\geq y)=\bar{\BE}[I_{|B_t|\geq y}]=\bar{\BE}[I_{|B_1|\geq\frac y{\sqrt{t}}}]
=\sup_{\theta\in\Theta}E_P[I_{|\int_0^1\theta_sdW_s|\geq \frac
y{\sqrt{t}}}].$$ Hence, for any $s\leq t$, we have for almost surely
$y\in\BR^+$, $V(|B_s|\geq y)\leq V(|B_t|\geq y).$\qquad$\square$
\begin{lemma}\label{le5}
If $\underline{\sigma}>0,$ then the paths of G-Brownian motion are
q.s. nowhere differentiable.
\end{lemma}
\textbf{Proof.} Suppose $B.(\omega)$ is differentiable at $s$, then
there exists $\delta>0,l\geq1$, such that
$|B_t(\omega)-B_s(\omega)|<l|t-s|$ for any $|t-s|<\delta.$  From the
definition of quadratic variation process of G-Brownian motion
\cite{p2006}, it follows that
$$\langle B\rangle_{s+\delta}(\omega)-\langle B\rangle_{s}(\omega)=\lim_{\mu(\pi_t^n)\to 0}
\sum_{i=0}^{n}|B_{t_{i+1}}(\omega)-B_{t_i}(\omega)|^2\leq
\lim_{\mu(\pi_t^n)\to 0}l^2\sum_{i=0}^{n}|t_{i+1}-t_i|^2=0.$$
Meanwhile Peng in \cite{p2010} has show that $\langle
B\rangle_{s+\delta}(\omega)-\langle
B\rangle_{s}>\underline{\sigma}\delta>0$ q.s.. We deduce therefore
that $V\{\omega:B.(\omega)$ is differentiable at $s\}=0$. So the
lemma holds.\qquad$\square$
\section{The invariance principle of G-Brownian motion}
In this section we will consider the invariance principle of
G-Brownian motion for the LIL under G-expectation. First let us give
some basic notations. Let $B(t)_{t\geq 0}$ be the G-Brownian motion,
$B_t\sim\CN(0,[t\GX,t\GS])$. Define
$$\zeta_n(t)=(2n\log\log n)^{-{1\over2}}B(nt),\ \ \ \forall t\in[0,1],n\geq3.$$
Let $C([0,1])$ be the banach space of continuous maps from $[0,1]$ to $\BR$ endowed with the supremum norm $\|\cdot \|$, using the enuclidean norm in $\BR$. $\zeta_n$ is then a random variable with values in $C([0,1])$. For any $\beta\in\BR^+$, define
$$K_{\beta}:=\{x(\cdot):x\in C([0,1]),x(0)=0,\int_0^1|\dot{x(t)}|^2dt\leq\beta^2\}.$$
\begin{theorem}\label{th1}
Let $C(\zeta_n)$ denotes the cluster of sequence $(\zeta_n)_{n=3}^{\infty}$, then \\
(I) The sequence $(\zeta_n)_{n\geq3}$ is relatively norm-compact q.s..\\
(II) $v\{C(\zeta_n)\subseteq K_{\bar{\sigma}}\}=1.$\\
(III) $v\{C(\zeta_n)\supseteq K_{\underline{\sigma}}\}=1.$ \\
(IV) $\forall \beta\in[\underline{\sigma},\bar{\sigma}], V\{C(\zeta_n)=K_{\beta}\}=1.$
\end{theorem}
\textbf{Proof.}
(I) and (II). For any $\e>0$, let
$$\KSE:=\{x(\cdot):\ x\in C([0,1]),d(x,\KS)\leq\e\}.$$
Moreover, let $\eta_n$ be the random variable in $ C([0,1])$ obtained by interpolating the points $\zeta_n({i\over m})$ at ${i\over m}\ (i=1,\cdots,m)$, where $m$
is a positive integer which will be decided in the later. For any $\e_1>0$, we have
\begin{eqnarray*}
&& V\{\zeta_n\notin\KSE\} \\
&&\leq V\{\eta_n\notin K_{\bar{\sigma}+\e_1}\}+V\{\eta_n\in K_{\bar{\sigma}+\e_1},\ \zeta_n\notin\KSE\} \\
&&=I_1+I_2.
\end{eqnarray*}
For any $\lambda>0$, by Chebyshev's inequality and Lemma \ref{le40},
we get
\begin{eqnarray*}
&&I_1=V\{\int_0^1|\dot{\eta_n(t)}|^2dt>(\bar{\sigma}+\e_1)^2\}\\
&&=V\left\{\frac{\sum_{i=1}^m [B(\frac{ni}m)-B(\frac{n(i-1)}m)]^2}{2n\log\log n/m}>(\bar{\sigma}+\e_1)^2\right\}\\
&&=V\left\{\frac{\sum_{i=1}^m [B(\frac{ni}m)-B(\frac{n(i-1)}m)]^2}{n\GS/m}>2(1+\frac{\e_1}{\bar{\sigma}})^2\log\log n\right\}\\
&&\leq\exp(-2\lambda(1+\frac{\e_1}{\bar{\sigma}})^2\log\log n) \bar{\BE}\left[\exp\left(\frac{\sum_{i=1}^m \lambda[B(\frac{ni}m)-B(\frac{n(i-1)}m)]^2}{n\GS/m}\right)\right]\\
&&=\exp(-2\lambda(1+\frac{\e_1}{\bar{\sigma}})^2\log\log n)
\prod_{i=1}^{m}\bar{\BE}\left[\exp\left(\frac{\lambda[B(\frac{ni}m)-B(\frac{n(i-1)}m)]^2}{n\GS/m}\right)\right]\\
&&=\exp(-2\lambda(1+\frac{\e_1}{\bar{\sigma}})^2\log\log n)
[\bar{\BE}\exp(\frac{\lambda B(\frac{n}m)^2}{n\GS/m})]^m
\end{eqnarray*}
For each $\e_1>0$, we choose $\lambda(\e_1)\in(0,\frac12)$ such that
$\beta_1:=2\lambda(1+\frac{\e_1}{\bar{\sigma}})^2>1$, therefore
$$C(\e_1):=\bar{\BE}\exp(\frac{\lambda B(\frac{n}m)^2}{n\GS/m})=\frac1{\sqrt{2\pi}}\int_{-\infty}^{\infty}\exp(\lambda y^2)\exp(-y^2/2)dy<\infty.$$
Hence
\begin{equation}\label{eq4}
I_1\leq C(\e_1)^m\exp(-\beta_1\log\log n).
\end{equation}
Meanwhile,
\begin{eqnarray*}
&&I_2=V\{\eta_n\in K_{\bar{\sigma}+\e_1},\ \zeta_n\notin\KSE\}\\
&&\leq V\{\eta_n\in K_{\bar{\sigma}+\e_1},\ \|\zeta_n-\frac{\bar{\sigma}}{\bar{\sigma}+\e_1}\eta_n\|\geq\e\}\\
&&=\sup_{P\in\CP}P\{\eta_n\in K_{\bar{\sigma}+\e_1},\
\|\zeta_n-\frac{\bar{\sigma}}{\bar{\sigma}+\e_1}\eta_n\|\geq\e\}
\end{eqnarray*}
Define the random variable $T$ by
\begin{displaymath}
T:=
\left\{\begin{array}{c}
\min\{t:\ t\in[0,1],|\zeta_n(t)-\frac{\bar{\sigma}}{\bar{\sigma}+\e_1}\eta_n(t)|\geq\e\},\ \ \textrm{if this set is nonempty}  ;\\2,\qquad\qquad\qquad\qquad
\qquad\qquad\ \ \ \textrm{otherwise},
\end{array} \right.
\end{displaymath}
and let $F_P$ be its distribution function under $P$, thus
\begin{eqnarray*}
&&I_2\leq\sup_{P\in\CP}\int_{0}^1P\{\eta_n\in K_{\bar{\sigma}+\e_1}|T=t\}dF_P(t)\\
&&=\sup_{P\in\CP}\int_{0}^1P\{\eta_n\in K_{\bar{\sigma}+\e_1},\
|\zeta_n(t)-\frac{\bar{\sigma}}{\bar{\sigma}+\e_1}\eta_n(t)|=\e|T=t\}dF_P(t).
\end{eqnarray*}
Let $i(t)$ denote the smallest integer $i$ with $i/m\geq t$, the statement $\eta_n\in K_{\bar{\sigma}+\e_1}$ implies
$$|\eta_n(\frac{i(t)}m)-\eta_n(t)|\leq\frac{\bar{\sigma}+\e_1}{\sqrt{m}}.$$
Together with $$|\zeta_n(t)-\frac{\bar{\sigma}}{\bar{\sigma}+\e_1}\eta_n(t)|=\e,$$
we have
\begin{eqnarray*}
&&|\zeta_n(\frac{i(t)}m)-\zeta_n(t)|\\
&&\geq|\eta_n(t)-\zeta_n(t)|-|\eta_n(\frac{i(t)}m)-\eta_n(t)|\\
&&\geq\frac{\bar{\sigma}+\e_1}{\bar{\sigma}}|\frac{\bar{\sigma}}{\bar{\sigma}+\e_1}\eta_n(t)-\zeta_n(t)-
\frac{-\e_1}{\bar{\sigma}+\e_1}\zeta_n(t)|-\frac{\bar{\sigma}+\e_1}{\sqrt{m}}\\
&&\geq(1+\frac{\e_1}{\bar{\sigma}})\e-\frac{\e_1|B(nt)|}{\bar{\sigma}\sqrt{2n\log\log n}}-\frac{\bar{\sigma}+\e_1}{\sqrt{m}}\\
&&\geq\e/2 \qquad q.s.
\end{eqnarray*}
where the last inequality is obtained by LIL of chen \cite{lil} and
choose $\e_1$ close to $0$ and $m$ be sufficiently large. For
$\lambda>0$, using lemma \ref{le3} it follows that
\begin{eqnarray*}
&&I_2\leq\sup_{P\in\CP}\int_{0}^1P\{|\zeta_n(\frac{i(t)}m)-\zeta_n(t)|\geq\e/2|T=t\}dF_P(t)\\
&&\leq \sup_{P\in\CP}V\{|\zeta_n({1\over m})|\geq\e/2\}\int_0^1dF_P(t)\\
&&\leq V\left\{\frac{|B({n\over m})|}{\sqrt{\frac{n\GS}m}}\geq\frac{\e\sqrt{2m\log\log n}}{2\bar{\sigma}}\right\}\\
&&\leq\exp(-\frac{\lambda\e^2m\log\log
n}{2\GS})\bar{\BE}[\exp(\frac{\lambda B({n\over
m})^2}{\frac{n\GS}m})].
\end{eqnarray*}
Choosing $\lambda(\e)\in(0,1/2)$ and $m$ such that $\beta_2:=\lambda\e^2m/\GS>1,$ as before there exist $D(\e)>0$ such that
\begin{equation}\label{eq5}
I_2\leq D(\e)\exp(-\beta_2\log\log n).
\end{equation}
From inequalities (\ref{eq4}) and (\ref{eq5})we have
$\beta:=\beta_1\wedge\beta_2>1$ and
$$V\{\zeta_n\notin\KSE\}\leq(C(\e_1)^m+D(\e))\exp(-\beta\log\log n).$$
If $n_k=\lfloor c^k \rfloor+1$, where $c>1,\lfloor c\rfloor$ is the
largest integer not greater than $c$, then
$$\sum_{k=1}^{\infty}V\{\zeta_{n_k}\notin\KSE\}\leq (C(\e_1)^m+D(\e))(\log c)^{-\beta}\sum_{k=1}^{\infty}k^{-\beta}<\infty.$$
By Borel-Cantelli lemma, we have
$$V\{\bigcap_{i=1}^{\infty}\bigcup_{k=i}^{\infty}(\zeta_{n_k}\notin\KSE)\}=0,$$
in other words,
$$v\{\bigcup_{i=1}^{\infty}\bigcap_{k=i}^{\infty}(\zeta_{n_k}\in\KSE)\}=1.$$
For $c$ sufficient close to 1 this implies that
$$v\{\bigcup_{i=1}^{\infty}\bigcap_{k=i}^{\infty}(\zeta_n\in K_{\bar{\sigma}}^{2\e})\}=1.$$
Hence, $v(C(\zeta_n)\subseteq K_{\bar{\sigma}})=1$ and for any $\e>0$, $(\zeta_n)_{n\geq3}$ exists a relatively compact $2\e$-net
q.s.. This shows that $(\zeta_n)_{n\geq3}$ is relatively compact q.s..
The proof of (I) and (II) is complete.
\par
(III). For any $0<\beta<\underline{\sigma},$ any given $x\in \KB$, let $m\geq 1$ be an integer, $\e_0>0$ and $A_n=\cap_{i=1}^mA_n^i$, where
$$A_n^i=\left\{\left|\zeta_n({i\over m})-\zeta_n({{i-1}\over m})-\left(x({i\over m})-x({{i-1}\over m})\right)\right|<\e_0\right\},\ \ \ i=1,\cdots,m.$$
By the stationary increments property of G-Brownian motion and Lemma
\ref{le40}, we have
\begin{eqnarray*}
v(A_n^i)&=&v\left\{\left|\frac{B(\frac nm)}{\sqrt{2n\log\log n}}-\left(x({i\over m})-x({{i-1}\over m})\right)\right|\leq\e_0\right\}\\
&=&v\left\{\left|\frac{\bar{B}(n)}{\sqrt{2n\log\log n}}-\left(x({i\over m})-x({{i-1}\over m})\right)\right|\leq\e_0\right\},
\end{eqnarray*}
where $\bar{B}(t)=B(\frac tm)\sim\CN(0,[\frac{\GX t}m,\frac{\GS
t}m])$. Let us choose $n_k=k^{k^{\alpha}}$ for $k\geq1$ where
$0<\alpha<{1\over 2m}$. Then we obtain that
\begin{eqnarray*}
v(A_{n_k}^i)&\geq& v\left\{\left|\frac{\BB(n_k)-\BB(n_{k-1})}{\sqrt{2n_k\log\log n_k}}-\left(x({i\over m})-x({{i-1}\over m})\right)\right|\leq\frac{\e_0}2\right\}\\
&&\cdot v\left\{\frac{|\BB(n_{k-1})|}{\sqrt{2n_{k-1}\log\log n_{k-1}}}\frac{\sqrt{2n_{k-1}\log\log n_{k-1}}}{\sqrt{2n_k\log\log n_k}}\leq\frac{\e_0}2\right\}
\end{eqnarray*}
For each $t>0$, let
$$N_k:=\lfloor(n_{k+1}-n_k)^2t^2/(2n_{k+1}\log\log n_{k+1})\rfloor,$$
$$l_k:=\lfloor2t^{-2}n_{k+1}\log\log n_{k+1}/(n_{k+1}-n_k)\rfloor,$$
$$r_k:=\sqrt{2n_{k+1}\log\log n_{k+1}}/(tl_k).$$
We have
\begin{eqnarray*}
&&v\left\{\left|\frac{\BB(n_k)-\BB(n_{k-1})}{\sqrt{2n_k\log\log n_k}}-\left(x({i\over m})-x({{i-1}\over m})\right)\right|\leq\frac{\e_0}2\right\}\\
&=&v\left\{x({i\over m})-x({{i-1}\over m})-\frac{\e_0}2\leq\frac{\BB(n_k-n_{k-1})}{\sqrt{2n_k\log\log n_k}}\leq x({i\over m})-x({{i-1}\over m})+\frac{\e_0}2\right\}\\
&\geq&v\left\{x({i\over m})-x({{i-1}\over m})-\frac{\e_0}4\leq\frac{\BB(N_{k-1}l_{k-1})}{\sqrt{2n_k\log\log n_k}}\leq x({i\over m})-x({{i-1}\over m})+\frac{\e_0}4\right\}\\
&\cdot&v\left\{-\frac{\e_0}4\leq\frac{\BB(n_k-n_{k-1})-\BB(N_{k-1}l_{k-1})}{\sqrt{2n_k\log\log n_k}}\leq \frac{\e_0}4\right\}\\
&\geq&v\left\{tl_{k-1}\left(x({i\over m})-x({{i-1}\over m})-\frac{\e_0}4\right)\leq\frac{\BB(N_{k-1}l_{k-1})}{r_{k-1}}\leq tl_{k-1}\left(x({i\over m})-x({{i-1}\over m})+\frac{\e_0}4\right)\right\}\\
&\cdot&v\left\{-\frac{\e_0}4\leq\frac{\BB(n_k-n_{k-1}-N_{k-1}l_{k-1})}{\sqrt{2n_k\log\log n_k}}\leq \frac{\e_0}4\right\}\\
&\geq&v\left\{\left(x({i\over m})-x({{i-1}\over m})\right)t-\frac{\e_0 t}4\leq\frac{\BB(N_{k-1})}{r_{k-1}}\leq \left(x({i\over m})-x({{i-1}\over m})\right)t+\frac{\e_0 t}4\right\}^{l_{k-1}}\\
&\cdot&v\left\{-\frac{\e_0}4\leq\frac{\BB(n_k-n_{k-1}-N_{k-1}l_{k-1})}{\sqrt{2n_k\log\log n_k}}\leq \frac{\e_0}4\right\}\\
&\geq&\left(\CE\left[\phi\left(\frac{\BB(N_{k-1})}{r_{k-1}}-(x({i\over m})-x({{i-1}\over m}))t\right)\right]\right)^{l_{k-1}}\\
&\cdot& v\left\{-\frac{\e_0}4\leq\frac{\BB(n_k-n_{k-1}-N_{k-1}l_{k-1})}{\sqrt{2n_k\log\log n_k}}\leq \frac{\e_0}4\right\},
\end{eqnarray*}
where $\phi(x)$ is a even function defined by
\begin{displaymath}
\phi(x):=
\left\{\begin{array}{c}
1-e^{|x|-\e_0 t/4},\ \ \ |x|\leq\e_0 t/4;\\0,\ \ \ \ \ \ \ \ \ \ \ \ \ \ \ \ |x|>\e_0 t/4.
\end{array} \right.
\end{displaymath}
Applying Lemma \ref{le2} and CLT \cite{p2009}, we have if
$k\to\infty$
\begin{eqnarray*}
&\log&\CE\left[\phi \left(\frac{\BB(N_{k-1})}{r_{k-1}}-(x({i\over m})-x({{i-1}\over m}))t\right)\right]\to\log\CE
\left[\phi\left(\BB(1)-(x({i\over m})-x({{i-1}\over m}))t\right)\right]\\
&\geq& -\frac{m(x({i\over m})-x({{i-1}\over m}))^2t^2}{2\GX}+\log\CE[\phi(\BB(1))].
\end{eqnarray*}
Thus as $k\to\infty$
\begin{eqnarray*}
&&\frac{n_k-n_{k-1}}{2n_k\log\log n_k}\cdot\log\CE\left[\phi\left(\frac{\BB(N_{k-1})}{r_{k-1}}-(x({i\over m})-x({{i-1}\over m}))t\right)\right]^{l_{k-1}}\\
&=&\frac{l_{k-1}(n_k-n_{k-1})}{2n_k\log\log n_k}\cdot \log\CE\left[\phi\left(\frac{\BB(N_{k-1})}{r_{k-1}}-(x({i\over m})-x({{i-1}\over m}))t\right)\right]\\
&\to&t^{-2}\log\CE[\phi(\BB(1)-(x({i\over m})-x({{i-1}\over m}))t)]\\
&\geq& -\frac{m(x({i\over m})-x({{i-1}\over m}))^2}{2\GX}+t^{-2}\log\CE[\phi(\BB(1))].
\end{eqnarray*}
Together with $\mathop {\underline {\lim } }\limits_{t \to \infty }t^{-2}\log\CE[\phi(\BB(1)-(x({i\over m})-x({{i-1}\over m}))t)]\geq  -\frac{m(x({i\over m})-x({{i-1}\over m}))^2}{2\GX},$
we have, for large enough $t$,
$$\lim_{k \to \infty }\frac{n_k-n_{k-1}}{2n_k\log\log n_k}\cdot \log\CE\left[\phi
\left(\frac{\BB(N_{k-1})}{r_{k-1}}-(x({i\over m})-x({{i-1}\over m}))t\right)\right]^{l_{k-1}}$$
\begin{equation}\label{eq2}
\geq-\frac{m(x({i\over m})-x({{i-1}\over m}))^2}{2\GX}-1.
\end{equation}
On the other hand, by Chebyshev's inequality, it follows that
$$V(|\frac{\BB(n_k-n_{k-1}-N_{k-1}l_{k-1}))}{\sqrt{2n_k\log\log n_k}}|>\frac{\e_0}4)
\leq \frac{8(n_k-n_{k-1}-N_{k-1}l_{k-1})\GS}{m\e_0^2n_k\log\log n_k}\to0,\ k\to\infty.$$
Therefore, as $k\to\infty,$
$$\frac{n_k-n_{k-1}}{2n_k\log\log n_k}\cdot\log v\left\{-\frac{\e_0}4\leq\frac{\BB(n_k-n_{k-1}-N_{k-1}l_{k-1}))}{\sqrt{2n_k\log\log n_k}}\leq \frac{\e_0}4\right\}$$
\begin{equation}\label{eq3}
=\frac{n_k-n_{k-1}}{2n_k\log\log n_k}\cdot\log\left(1-V(|\frac{\BB(n_k-n_{k-1}-N_{k-1}l_{k-1}))}{\sqrt{2n_k\log\log n_k}}|>\frac{\e_0}4)\right)\to0.
\end{equation}
So, from (\ref{eq2}) and (\ref{eq3}) we have
\begin{eqnarray*}
&&\mathop {\underline {\lim } }\limits_{k \to \infty }\frac{n_k-n_{k-1}}{2n_k\log\log n_k}\cdot\log v\left\{\left|\frac{\BB(n_k)-\BB(n_{k-1})}{\sqrt{2n_k\log\log n_k}}-(x({i\over m})-x({{i-1}\over m}))\right|\leq\frac{\e_0}2\right\}\\
&\geq&-\frac{m(x({i\over m})-x({{i-1}\over m}))^2}{2\GX}-1.
\end{eqnarray*}
Since $x\in\KB,$ we have
$$|x({i\over m})-x({{i-1}\over m})|\leq\frac{\beta}{\sqrt{m}}<\frac{\underline{\sigma}}{\sqrt{m}},\ \ \ \forall i=1,\cdots,m.$$
So there exist $\delta>0$ such that $d:=\delta+\max_{i\leq
m}(\frac{x({i\over m})-x({{i-1}\over
m})}{\frac{\underline{\sigma}}{\sqrt m}})^2+1<2.$ Thus, there exist
$k_0$ such that $\forall k\geq k_0$,
\begin{eqnarray*}
&&v\left\{|\frac{\BB(n_k)-\BB(n_{k-1})}{\sqrt{2n_k\log\log n_k}}-(x({i\over m})-x({{i-1}\over m}))|\leq\frac{\e_0}2\right\}\\
&&\geq \exp({-dn_k\log\log n_k/(n_k-n_{k-1})}).
\end{eqnarray*}
Meanwhile from the LIL of chen \cite{lil} we have $v(\mathop {\overline {\lim } }\limits_{n \to \infty }\frac{\BB(n)}{\sqrt{2n\log\log n}}\leq\frac{\bar{\sigma}}{\sqrt{m}})=1$ and $\frac{\sqrt{2n_{k-1}\log\log n_{k-1}}}{\sqrt{2n_k\log\log n_k}}\rightarrow0$ as $k\rightarrow\infty$, hence there exist $k_1$ such that $\forall k\geq k_1$,
\begin{eqnarray*}
 v\left\{\frac{|\BB(n_{k-1})|}{\sqrt{2n_{k-1}\log\log n_{k-1}}}\frac{\sqrt{2n_{k-1}\log\log n_{k-1}}}{\sqrt{2n_k\log\log n_k}}<\frac{\e_0}2\right\}\geq{1\over2}.
\end{eqnarray*}
By applying lemma \ref{le40}, we get for any $k>(k_0\vee k_1)$
\begin{eqnarray*}
v(A_{n_k})&=&\prod_{i=1}^m v(A_{n_k}^i)\\
&\geq&{1\over{2^m}}\exp({-mdn_k\log\log n_k/(n_k-n_{k-1})})\\
&\geq&{1\over{2^m}}\exp({-2mn_k\log\log n_k/(n_k-n_{k-1})})\\
&\sim&{1\over{2^m}}\exp(-2m\log\log n_k)\\
&\sim&\frac1{k^{2m\alpha}(2\log k)^{2m}}.
\end{eqnarray*}
Thus $\sum_{k=1}^{\infty}v(A_{n_k})=\infty$ for $2m\alpha<1$, using
the Borel-Cantelli lemma, we get infinitely many events $A_{n_k}$
happen q.s..
\par
Next will show that for any given $0<s<t<1$, quasi-surely there is infinitely
$k$ such that
\begin{equation}\label{eq1}
\zeta_{n_k}(t)-\zeta_{n_k}(s)\leq\bar{\sigma}\sqrt{t-s}+\e_0.
\end{equation}
In fact, by the LIL under capacity, we have
$$v(\mathop {\overline {\lim } }\limits_{n \to \infty }\frac{|\tilde{B}(n)|}{\sqrt{2n\log\log n}}\leq\bar{\sigma}\sqrt{t-s})=1,$$
where $\tilde{B}(r)=B((t-s)r)\sim\CN(0,[(t-s)\GX r,(t-s)\GS r])$.
So, we have
\begin{eqnarray*}
&&v\{\bigcap_{i=1}^{\infty}\bigcup_{k=i}^{\infty}(|\zeta_{n_k}(t)-\zeta_{n_k}(s)|\leq\bar{\sigma}\sqrt{t-s}+\e_0)\}\\
&&=v\{\bigcap_{i=1}^{\infty}\bigcup_{k=i}^{\infty}(\frac{|\tilde{B}(n_k)|}{\sqrt{2n_k\log\log n_k}}\leq\bar{\sigma}\sqrt{t-s}+\e_0)\}=1.
\end{eqnarray*}
For any $\varepsilon>0$ we set $m>(\frac{4\bar{\sigma}}{\e})^2$ and choose $\e_0$ such that $\e_0<\frac{\e}{2m}$ after fixed $m$. Now we consider the $n_k$ making $A_{n_k}$ happens and satisfying inequality (\ref{eq1}),
\begin{eqnarray*}
&&\|\zeta_{n_k}-x\|=\sup_{t\in[0,1]}|\zeta_{n_k}(t)-x(t)|\\
&&=\sup_{t\in[0,1]}\left|\zeta_{n_k}(t)-\zeta_{n_k}(\frac{\lfloor mt\rfloor}m)+x(\frac{\lfloor mt\rfloor}m)-x(t)\right.\\
&&\left.+\sum_{i=1}^{\lfloor mt\rfloor}\left[\zeta_{n_k}(\frac im)-\zeta_{n_k}(\frac{i-1}m)-\left(x(\frac im)-x(\frac{i-1}m)\right)\right]\right|\\
&&\leq\frac{\bar{\sigma}}{\sqrt m}+\frac{\beta}{\sqrt m}+m\e_0\\
&&\leq\frac{2\bar{\sigma}}{\sqrt m}+m\e_0\\
&&\leq\e.
\end{eqnarray*}
We conclude that $$v\{x\in C(\zeta_n)\}=1,$$ and thus
$$v\{\KB\subseteq C(\zeta_n)\}=1,$$ for $\KB$ having
countable dense set. So $$v\{\bigcup_{n=1}^{\infty}
K_{\underline{\sigma}-{1\over n}}\subseteq C(\zeta_n)\}=1.$$ Since
$C(\zeta_n))$ is a closed set, we get
$$v\{K_{\underline{\sigma}}\subseteq C(\zeta_n)\}=1.$$
The proof of (III) is complete.
\par
(IV). For any $\beta\in[\underline{\sigma},\bar{\sigma}]$, there exist $P_{\beta}\in\CP$ such that $B(t)/\beta$ be a classical Brownian motion under
$P_{\beta}$, so by the stranssen's invariance principle $P_{\beta}(C(B_n)=K_\beta)=1$. Therefore, $V(C(B_n)=K_\beta)=1.$
\par
The proof of Theorem \ref{th1} is complete.\qquad$\square$
\par
Notice that the discreteness of $n$ is inessential for the previous
considerations. More precisely, the following corollary holds true.
\begin{corol}\label{co1}
 If $u>e$ is real and we put $\zeta_u(t)=(2u\log\log u)^{-{1\over2}}B(ut),\ t\in[0,1]$, then we have\\
(I) The sequence $(\zeta_u)_{u>e}$ is relatively norm-compact q.s..\\
(II) $v\{C(\zeta_u)\subseteq K_{\bar{\sigma}}\}=1.$\\
(III) $v\{C(\zeta_u)\supseteq K_{\underline{\sigma}}\}=1.$ \\
(IV) $\forall \beta\in[\underline{\sigma},\bar{\sigma}], V\{C(\zeta_u)=K_{\beta}\}=1.$
\end{corol}
\begin{corol}\label{co2}
If $\varphi$ is a continuous map from $C[0,1]$ to some Hausdorff space $H$, then we have\\
(I) The sequence $(\varphi(\zeta_n))_{n\geq3}$ is relatively norm-compact q.s..\\
(II) $v\{C(\varphi(\zeta_n))\subseteq \varphi(K_{\bar{\sigma}})\}=1.$\\
(III) $v\{C(\varphi(\zeta_n))\supseteq \varphi(K_{\underline{\sigma}})\}=1.$ \\
(IV) $\forall \beta\in[\underline{\sigma},\bar{\sigma}], V\{C(\varphi(\zeta_n))=\varphi(K_{\beta})\}=1.$ \\
We substitute $(\zeta_n)_{n\geq3}$ with $(\zeta_u)_{u>e}$, the conclusion also holds.
\end{corol}
\section{Some applications and comments}
In this section, we give some applications and comments, which can be obtained by our invariance principle and the arguments of Strassen \cite{strassen}.
\begin{example}
Let $f(\cdot)$ be any Riemann integrable real function on $[0,1]$,
$$F(t)=\int_t^1 f(s)ds,\ t\in[0,1].$$
Then, $$v\left\{\underline{\sigma}\left(\int_0^1F^2(t)
dt\right)^{1/2}\leq \mathop {\overline {\lim } }\limits_{n \to
\infty }(2n^3\log\log n)^{-1/2}\sum_{i=1}^n f({i\over n})B_i\leq
\bar{\sigma}\left(\int_0^1F^2(t) dt\right)^{1/2}\right\}=1.$$ In
particular, for any $\alpha>-1$, putting $f(t)=t^{\alpha}$, we have
$$v\left\{\frac{\underline{\sigma}}{\sqrt{(\alpha+3/2)(\alpha+2)}}\leq \mathop {\overline {\lim } }\limits_{n \to \infty }(2n^{2\alpha+3}\log\log
n)^{-1/2}\sum_{i=1}^n i^{\alpha}B_i\leq
\frac{\bar{\sigma}}{\sqrt{(\alpha+3/2)(\alpha+2)}}\right\}=1.$$
\end{example}
\begin{example}
Let $a\geq1$ be real, then we have
$$v\left\{\frac{2(a+2)^{(a/2)-1}}{a^{a/2}\left(\int_0^1\frac{dt}{\underline{\sigma}\sqrt{1-t^a}}\right)^a}
\leq \mathop {\overline {\lim } }\limits_{n \to \infty
}n^{-1-(a/2)}(2\log\log n)^{-1/2}\sum_{i=1}^n|B_i|^a \leq
\frac{2(a+2)^{(a/2)-1}}{a^{a/2}\left(\int_0^1\frac{dt}{\bar{\sigma}\sqrt{1-t^a}}\right)^a}
\right\}=1.$$ In particular, a=1,2, we have
$$v\left\{\frac{\underline{\sigma}}{\sqrt 3}\leq \mathop {\overline {\lim } }\limits_{n \to \infty
}n^{-3/2}(2\log\log n)^{-1/2}\sum_{i=1}^n|B_i| \leq
\frac{\bar{\sigma}}{\sqrt 3}\right\}=1,$$
$$v\left\{\frac{4\GX}{\pi^2}\leq \mathop {\overline {\lim } }\limits_{n \to \infty
}n^{-2}(2\log\log n)^{-1}\sum_{i=1}^n|B_i|^2 \leq
\frac{4\GS}{\pi^2}\right\}=1.$$
\end{example}

\end{document}